\documentstyle[11pt,twoside]{article}
\include{graphicx}
\title{A uniform asymptotic expansion for the incomplete gamma functions revisited}
\author{\sc R. B.\ Paris \\
{\em Division of Computing and Mathematics}, \\
{\em University of Abertay Dundee, Dundee DD1 1HG, UK}
}
\begin{document}
\def\f#1#2{\mbox{${\textstyle \frac{#1}{#2}}$}}
\def\dfrac#1#2{\displaystyle{\frac{#1}{#2}}}
\def\boldal{\mbox{\boldmath $\alpha$}}
{\newcommand{\Sgoth}{S\;\!\!\!\!\!/}
\newcommand{\bee}{\begin{equation}}
\newcommand{\ee}{\end{equation}}
\newcommand{\lam}{\lambda}
\newcommand{\ka}{\kappa}
\newcommand{\al}{\alpha}
\newcommand{\fr}{\frac{1}{2}}
\newcommand{\fs}{\f{1}{2}}
\newcommand{\g}{\Gamma}
\newcommand{\br}{\biggr}
\newcommand{\bl}{\biggl}
\newcommand{\ra}{\rightarrow}
\newcommand{\mbint}{\frac{1}{2\pi i}\int_{c-\infty i}^{c+\infty i}}
\newcommand{\mbcint}{\frac{1}{2\pi i}\int_C}
\newcommand{\mboint}{\frac{1}{2\pi i}\int_{-\infty i}^{\infty i}}
\newcommand{\gtwid}{\raisebox{-.8ex}{\mbox{$\stackrel{\textstyle >}{\sim}$}}}
\newcommand{\ltwid}{\raisebox{-.8ex}{\mbox{$\stackrel{\textstyle <}{\sim}$}}}
\renewcommand{\topfraction}{0.9}
\renewcommand{\bottomfraction}{0.9}
\renewcommand{\textfraction}{0.05}
\newcommand{\mcol}{\multicolumn}
\date{}
\maketitle
\pagestyle{myheadings}
\markboth{\hfill \sc R. B.\ Paris  \hfill}
{\hfill \sc  Incomplete gamma function expansions\hfill}
\begin{abstract}
A new uniform asymptotic expansion for the incomplete gamma function $\Gamma(a,z)$ valid for large values of $z$ was 
given by the author in {\it J. Comput. Appl. Math.} {\bf 148} (2002) 323--339. This expansion contains a complementary error function of an argument measuring transition across the point $z=a$, with easily computable coefficients that do not involve a removable singularity 
in the neighbourhood of this point. In this note we correct a misprint in the listing of certain coefficients in this expansion and discuss in more detail the situation corresponding to $\Gamma(a,a)$.
\vspace{0.4cm}

\vspace{0.3cm}

\end{abstract}

\vspace{0.3cm}

\noindent $\,$\hrulefill $\,$

\vspace{0.2cm}

\begin{center}
{\bf 1. \  Introduction}
\end{center}
\setcounter{section}{1}
\setcounter{equation}{0}
\renewcommand{\theequation}{\arabic{section}.\arabic{equation}}
The incomplete gamma functions $\gamma(a,z)$ and $\g(a,z)$ (together with their normalised counterparts $P(a,z)$ and $Q(a,z)$) are defined by the integrals
\begin{eqnarray}
P(a,z)\!\!&:=&\!\!\frac{\gamma(a,z)}{\g(a)}=\frac{1}{\g(a)} \int_0^z t^{a-1}e^{-t} dt\qquad (\Re (a)>0),\nonumber\\
Q(a,z)\!\!&:=&\!\!\frac{\Gamma(a,z)}{\g(a)}=\frac{1}{\g(a)} \int_z^\infty t^{a-1}e^{-t} dt \qquad (|\arg\,z|<\pi),\label{e11}
\end{eqnarray}
where $a$ and $z$ are complex variables and the integration paths do not cross the negative real axis. The normalised functions $P(a,z)$ and $Q(a,z)$ satisfy the identity
\bee\label{e11a}
P(a,z)+Q(a,z)=1.
\ee
A review of the different types of expansion for $\g(a,z)$ that exist in the literature is presented in \cite{P}.
Earlier expansions not valid in the neighbourhood of the transition point $z=a$ are associated with the names of Mahler
(1930) \cite{M} and Tricomi (1950) \cite{T}, although Tricomi also gave an expansion of a more uniform character for $Q(a+1,a+x\sqrt{2a})$. 

A significant advance in the understanding of the uniform asymptotic structure of the incomplete gamma function was made by Temme \cite{Tem}, where it was established that
\bee\label{e12}
Q(a,z)\sim \frac{1}{2} \mbox{erfc} \bl(\eta \sqrt{\fs a}\br)+\frac{e^{-\fr a\eta^2}}{\sqrt{2\pi a}}\sum_{k=0}^\infty
{\cal C}_k(\eta)\,a^{-k}
\ee
as $a\ra\infty$. The auxiliary variable $\eta$ measures transition through $z=a$ and is defined by
$\eta=\{2(\lambda-1-\log\,\lambda)\}^{1/2}$, $\lambda=z/a$, with the branch being chosen so that $\eta(\lambda)$ is analytic near $\lambda=1$ ($z=a$) and $\eta\sim\lambda-1$ as $\lambda\ra 1$. The coefficients ${\cal C}_k(\eta)$ are defined recursively by
\[{\cal C}_0(\eta)=\frac{1}{\lambda-1}-\frac{1}{\eta},\qquad {\cal C}_k(\eta)=\frac{1}{\eta}\,\frac{d}{d\eta}\,{\cal C}_{k-1}(\eta)+\frac{\gamma_k}{\lambda-1} \quad (k\geq 1)\]
where $\gamma_k$ are the Stirling coefficients; see (\ref{e30}). The expansion (\ref{e12}) holds uniformly in $|z|\in [0,\infty)$ for a wide domain of complex $a$ and $z$; see \cite{Tem}. The above form captures the change in asymptotic structure through the point $z=a$ (where $\eta=0$) as well as describing the large-$a$ asymptotics, but suffers from the inconvenience
of the coefficients ${\cal C}_k(\eta)$ possessing a removable singularity at $\eta=0$. This last fact results in their evaluation in the neighbourhood of the transition point being difficult.

A new uniform asymptotic expansion for the incomplete gamma functions valid for large values of $z$ was given in \cite{P}. For $\g(a,z)$ this takes the form
\bee\label{e13}
\g(a,z)\sim z^{a-\fr}e^{-z}\bl\{\sqrt{\frac{\pi}{2}}\,e^{\chi^2/2} \mbox{erfc} \bl(\frac{\chi}{\sqrt{2}}\br) \sum_{k=0}^\infty\frac{A_k(\chi)}{z^{k/2}}-\sum_{k=1}^\infty\frac{B_k(\chi)}{z^{k/2}}\br\}
\ee
for $\Re (z-a)\geq 0$, where the variable $\chi=(z-a)/\sqrt{z}$, together with an analogous result for $\gamma(a,z)$.
This expansion also contains a complementary error function of an argument measuring transition through the point $z=a$, with easily computable coefficients $A_k(\chi)$ and $B_k(\chi)$ that do not involve a removable singularity at $z=a$.
An expansion of a similar nature for $\g(a,z)$, but for $a\ra\infty$, is described in Dingle's book \cite[p.~249]{D}; see also \cite[Section 4]{P}.

The procedure employed in \cite{D, P} consists of factorising the exponential factors appearing in the integrands of appropriate integral representations for $\g(a,z)$ and $\gamma(a,z)$ into an exponential factor containing only the linear and quadratic terms and another factor that is expanded in ascending powers of the integration variable.
Although these expansions are not valid in as large a domain of complex $a$ and $z$ as that in (\ref{e12}), the coefficients have the advantage of being more straightforward to compute near $z=a$.

In this note we correct an error in the presentation of the coefficients $A_k(\chi)$ and $B_k(\chi)$ in \cite{P}
and discuss their computation away from the neighbourhood of $z=a$. We also examine in more detail the expansion (\ref{e13}) when $z=a$.
\vspace{0.6cm}

\begin{center}
{\bf 2. \ The expansion for $\g(a,z)$}
\end{center}
\setcounter{section}{2}
\setcounter{equation}{0}
\renewcommand{\theequation}{\arabic{section}.\arabic{equation}}
For completeness in presentation we give an outline of the arguments described in \cite[Section 2]{P} omitting technical details.
We first consider the expansion for $\g(a,z)$ for large $z$ where, with the change of variable $t=ze^\tau$ in (\ref{e11}), we find
\bee\label{e21}
\g(a,z)=z^a\int_0^\infty e^{a\tau-ze^\tau} d\tau = z^ae^{-z}\int_0^\infty e^{-(z-a)\tau-\fr z\tau^2} H(\tau;z)\,d\tau
\ee
for $|\arg\,z|<\fs\pi$. In (\ref{e21}) we have written the exponential $e^{a\tau-ze^\tau}$ in factored form with
\[H(\tau;z)=e^{-zh(\tau)},\qquad h(\tau):=e^\tau-1-\tau-\fs\tau^2.\]

We expand the factor $H(\tau;z)$ as a finite Taylor series
\bee\label{e22}
H(\tau;z)=\sum_{k=0}^{n-1}\frac{c_k(z)}{k!}\,\tau^k+r_n(\tau;z)\qquad (n=1, 2, \ldots),
\ee
where $r_n(\tau;z)$ denotes the remainder term. The coefficients $c_k(z)$ satisfy the recurrence relation
\[c_{k+1}(z)=-z\sum_{j=2}^k\bl(\!\!\begin{array}{c}k\\j\end{array}\!\!\br)\,c_{k-j}(z)\qquad (k\geq 2),\]
with $c_0(z)=1$ and $c_1(z)=c_2(z)=0$. The $c_k(z)$ are polynomials in $z$ of degree $\lfloor k/3\rfloor$ 
given by
\bee\label{e23}
c_k(z)=\sum_{j=0}^{\lfloor k/3\rfloor} (-)^j S_3(k,j)\,z^j,\qquad S_3(k,0)=\delta_{0k},
\ee
where the coefficients\footnote{In \cite{P}, the coefficients $S_3(k,j)$ were denoted by $\alpha_j(k)$.} $S_3(k,j)$
are the 3-associated Stirling numbers of the second kind
defined by \cite[p.~222]{Com}
\[\exp \bl[u\bl(\frac{t^3}{3!}+\frac{t^4}{4!}+\cdots \br)\br]=\sum_{k=0}^\infty \sum_{j=0}^{\lfloor k/3 \rfloor} S_3(k,j) \,\frac{u^jt^k}{k!}\]
and $\delta_{0k}$ is the Kronecker symbol.
In Table 1 we show the values of $S_3(k,j)$ for $1\leq j\leq 6$ and $3\leq k\leq 20$.

\begin{table}[t]
\caption{\footnotesize{Values of the coefficients $S_3(k,j)$ for $1\leq j\leq 6$ and $3\leq k\leq 20$. }}
\begin{center}
\begin{tabular}{|r|llllll|}
\hline
&&&&&&\\[-0.3cm]
\mcol{1}{|c|}{$k\backslash j$} & \mcol{1}{c}{1} & \mcol{1}{c}{2} & \mcol{1}{c}{3} & \mcol{1}{c}{4} & \mcol{1}{c}{5} & \mcol{1}{c|}{6}\\
\hline
&&&&&&\\[-0.3cm]
3  & 1 &&&&&\\
4  & 1 &&&&&\\
5  & 1 &&&&&\\
6  & 1 & 10 &&&&\\
7  & 1 & 35 &&&&\\
8  & 1 & 91 &&&&\\
9  & 1 & 210 & 280 &&&\\
10 & 1 & 456 & 2100 &&&\\
11 & 1 & 957 & 10395 &&&\\
12 & 1 & 1969 & 42735 & 15400 &&\\
13 & 1 & 4004 & 158301 & 200200 &&\\
14 & 1 & 8086 & 549549 & 1611610 &&\\
15 & 1 & 16263 & 1827826 & 10335325 & 1401400 &\\
16 & 1 & 32631 & 5903898 & 57962905 & 28028000 &\\
17 & 1 & 65382 & 18682014 & 297797500 & 333533200 &\\
18 & 1 & 130900 & 58257810 & 1439774336 & 3073270200 & 190590400\\
19 & 1 & 261953 & 179765973 & 6662393738 & 24234675465 & 5431826400\\
20 & 1 & 524077 & 550478241 & 29844199346 & 172096749825 & 89625135600\\
[.2cm]\hline
\end{tabular}
\end{center}
\end{table}
Substitution of the expansion (\ref{e22}) into the integral (\ref{e21}) then yields
\bee\label{e24}
\g(a,z)=z^ae^{-z}\bl\{\sum_{k=0}^{n-1}\frac{c_k(z)}{k!} \int_0^\infty \tau^k e^{-(z-a)\tau-\fr z\tau^2}\,d\tau+R_n(a,z)\br\},
\ee
where
\[R_n(a,z)=\int_0^\infty e^{-(z-a)\tau-\fr z\tau^2} r_n(\tau;z)\,d\tau.\]
The integrals in (\ref{e24}) can be evaluated in terms of the parabolic cylinder function $D_\nu(z)$ as
\[\int_0^\infty \tau^k e^{-(z-a)\tau-\fr z\tau^2}\,d\tau=z^{-(k+1)/2} k!\,d_k(\chi),\qquad d_k(\chi):=e^{\chi^2/4} D_{-k-1}(\chi),\]
\[ \chi:=\frac{z-a}{\sqrt{z}},\]
whence we obtain the result
\[
\g(a,z)=z^{a}e^{-z}\bl\{\sum_{k=0}^{n-1}\frac{c_k(z)}{z^{(k+1)/2}} \,d_k(\chi)+R_n(a,z)\br\}.
\]
We now choose $n=3m+4$, $m=0, 1, 2 \ldots\ $. In \cite{P}, it is established that $R_n(a,z)=O(z^{-m/2-1})$ as $z\ra\infty$ in $|\arg\,z|<\fs\pi$. Thus we obtain\footnote{The terms in the sum appearing in (\ref{e26}) corresponding to $3m+1\leq k\leq 3m+3$ are O$(z^{-m/2-1})$ and so have been included in the order term.} 
\bee\label{e26}
\g(a,z)=z^{a-\fr}e^{-z}\bl\{\sum_{k=0}^{3m}\frac{c_k(z)}{z^{k/2}} \,d_k(\chi)+O(z^{-(m+1)/2})\br\},\qquad \Re (z-a)\geq 0
\ee
as $|z|\ra\infty$ in $|\arg\,z|<\fs\pi$. 

The functions $d_k(\chi)$ involve the parabolic cylinder functions of negative integer order and so can be expressed by recurrence in terms of the complementary error function, since
\bee\label{e26a}
d_0(\chi)=\sqrt{\frac{\pi}{2}}\, e^{\chi^2/2}\,\mbox{erfc} \bl(\frac{\chi}{\surd 2}\br).
\ee
From the recurrence relation $D_{\nu+1}(z)-zD_\nu(z)+\nu D_{\nu-1}(z)=0$ satisfied by the parabolic cylinder function, we obtain
\bee\label{e27}
d_{k+1}(\chi)=\frac{1}{k+1} \{d_{k-1}(\chi)-\chi\,d_k(\chi)\} \qquad (k=0, 1, 2 \ldots  ),
\ee
where $d_{-1}(\chi)=e^{\chi^2/4} D_0(\chi)=1$. This enables us to express $d_k(\chi)$ in the form
\bee\label{e27a}
d_k(\chi)=p_k(\chi) d_0(\chi)-q_k(\chi),
\ee
where $p_k(\chi)$, $q_k(\chi)$ are polynomials in $\chi$ of degree $k$ and $k-1$, respectively. The functions 
$p_k(\chi)$ and $q_k(\chi)$ satisfy the recurrence (\ref{e27}) with the starting values $(p_{-1}, p_0)=(0, 1)$ and $(q_{-1}, q_0)=(-1,0)$. The first few values are 
\[p_0(\chi)=1,\ \ p_1(\chi)=-\chi,\ \ p_2(\chi)=\fs+\fs\chi^2,\ \ p_3(\chi)=-\fs\chi-\f{1}{6}\chi^3,\ \ p_4(\chi)=\f{1}{8}+\f{1}{4}\chi^2+\f{1}{24}\chi^4,\]
\[q_0(\chi)=0,\ \ q_1(\chi)=-1,\ \ q_2(\chi)=\fs\chi,\ \ q_3(\chi)=-\f{1}{3}-\f{1}{6}\chi^2,\ \ q_4(\chi)=\f{5}{24}\chi^2+\f{1}{24}\chi^3,\]
and so on.

We substitute (\ref{e23}) into the expansion (\ref{e26}) and collect terms involving like powers of $z$ together to
yield
\[
\g(a,z)=z^{a-\fr}e^{-z}\bl\{\sum_{\ell=0}^{3m}  \frac{d_\ell(\chi)}{z^{\ell/2}}\sum_{j=0}^{\lfloor \ell/3\rfloor}(-)^j S_3(\ell,j) z^j+O(z^{-(m+1)/2})\br\}\]
\[= z^{a-\fr}e^{-z}
\bl\{\sum_{k=0}^{m}\frac{C_k(\chi)}{z^{k/2}}+O(z^{-(m+1)/2})\br\},\hspace{1.8cm}\]
where some terms have been incorporated into the order term. The coefficients $C_k(\chi)$ are defined by
\bee\label{e28}
C_k(\chi)=\sum_{j=0}^k (-)^j S_3(k+2j,j)\,d_{k+2j}(\chi)
= A_k(\chi) d_0(\chi)-B_k(\chi),
\ee
where
\bee\label{e29}
\left\{\!\!\begin{array}{l}A_k(\chi)\\B_k(\chi)\end{array}\!\!\right\}=\sum_{j=0}^k (-)^j S_3(k+2j,j) \left\{\!\!\begin{array}{r}p_{k+2j}(\chi)\\q_{k+2j}(\chi)\end{array}\!\!\right\}.
\ee
The coefficients $A_k(\chi)$ and $B_k(\chi)$ are polynomials in $\chi$ of degree $3k$ and $3k-1$, respectively, with $B_0(\chi)=0$; their expressions for $0\leq k\leq 4$ are presented in Table 2. In \cite[Table 2]{P} the expressions are given for $0\leq k\leq 8$. However, the table headings in \cite{P} are incorrect: the headings $(-)^k A_k(\chi)$ and $(-)^{k+1}B_k(\chi)$ should read $A_k(\chi)$ and $-B_k(\chi)$, respectively. Note also that in this paper we have reversed the sign of $B_k(\chi)$ for convenience in presentation. 

\begin{table}[th]
\caption{\footnotesize{Values of the coefficients $A_k(x)$ and $B_k(x)$ for $0\leq k\leq 5$. }}
\begin{center}
\begin{tabular}{|r|l|}
\hline
&\\[-0.3cm]
\mcol{1}{|c|}{$k$} & \mcol{1}{c|}{$A_k(x)$}\\
\hline
&\\[-0.3cm]
0 & 1\\
[0.15cm]
1 & $\fs x+\f{1}{6}x^3$\\
[0.15cm]
2 & $\f{1}{12}+\f{3}{8}x^2+\f{1}{6}x^4+\f{1}{72}x^6$\\
[0.1cm]
3 & $\f{1}{8}x+\f{47}{144}x^3+\f{37}{240}x^5+\f{1}{48}x^7+\f{1}{1296}x^9$\\
[0.15cm]
4 & $\f{1}{288}+\f{5}{32}x^2+\f{347}{1152}x^4+\f{617}{4320}x^6+\f{23}{960}x^8+\f{1}{648}x^{10}+\f{1}{31104}x^{12}$\\
[0.15cm]
5 & $\f{5}{576}x+\f{79}{432}x^3+\f{367}{1280}x^5+\f{32353}{241920}x^7+\f{785}{31104}x^9+\f{37}{17280}x^{11}+\f{5}{62208}x^{13}+\f{1}{933120}x^{15}$\\
[.2cm]\hline
&\\[-0.3cm]
\mcol{1}{|c|}{$k$} & \mcol{1}{c|}{$B_k(x)$}\\
\hline
&\\[-0.3cm]
0 & 0\\
[0.15cm]
1 & $\f{1}{3}+\f{1}{6}x^2$\\
[0.15cm]
2 & $\f{1}{4}x+\f{11}{72}x^3+\f{1}{72}x^5$\\
[0.15cm]
3 & $\f{4}{135}+\f{241}{1080}x^2+\f{293}{2160}x^4+\f{13}{648}x^6+\f{1}{1296}x^8$\\
[0.15cm]
4 & $\f{241}{4320}x+\f{341}{1620}x^3+\f{6377}{51840}x^5+\f{389}{17280}x^7+\f{47}{31104}x^9+\f{1}{31104}x^{11}$\\
[0.15cm]
5 & $-\f{8}{2835}+\f{14297}{181440}x^2+\f{7403}{36288}x^4+\f{9179}{80640}x^6+\f{403}{17280}x^{8}+\f{107}{51840}x^{10}+\f{37}{466560}x^{12}+\f{1}{933120}x^{14}$\\
[0.2cm]
\hline
\end{tabular}
\end{center}
\end{table}

Then we obtain the following theorem:
\newtheorem{theorem}{Theorem}
\begin{theorem}$\!\!\!.$ Let $m=0, 1, 2, \ldots\ $ and $\chi=(z-a)/\sqrt{z}$. Then, when $\Re (z-a)\geq 0$,  we have the expansion
\bee\label{e210}
\g(a,z)=z^{a-\fr}e^{-z}\bl\{d_0(\chi) \sum_{k=0}^m \frac{A_k(\chi)}{ z^{k/2}}-\sum_{k=1}^m \frac{B_k(\chi)}{ z^{k/2}} +O(z^{-(m+1)/2})\br\}
\ee
as $z\ra\infty$ in the sector $|\arg\,z|<\fs\pi$. The function $d_0(\chi)$ is given in (\ref{e26a}) and the coefficients
$A_k(\chi)$, $B_k(\chi)$ are defined in (\ref{e29}).
\end{theorem}
\vspace{0.2cm}

A similar treatment can be applied to the function $\gamma(a,z)$. Making the change of variable $t=ze^{-\tau}$ in (\ref{e11}) we obtain
\[\gamma(a,z)=z^a\int_0^\infty e^{-a\tau-ze^{-\tau}} d\tau=z^ae^{-z}\int_0^\infty e^{-(a-z)\tau-\fr z\tau^2} H(-\tau;z)\,d\tau\]
provided $\Re (a)>0$.
Substitution of the expansion for $H(-\tau;z)$ in this last integral then produces (see \cite{P})
\[
\gamma(a,z)=z^{a-\fr}e^{-z}\bl\{\sum_{k=0}^{3m}\frac{(-)^kc_k(z)}{z^{k/2}} \,d_k(-\chi)+O(z^{-\fr m-1})\br\},\qquad \Re (z-a)\leq 0
\]
as $|z|\ra\infty$ in $|\arg\,z|<\fs\pi$. This leads to the expansion given by \cite{P}
\begin{theorem}$\!\!\!.$ Let $m=0, 1, 2, \ldots\ $ and $\chi=(z-a)/\sqrt{z}$. Then, when $\Re (z-a)\leq 0$,  we have the expansion
\bee\label{e211}
\gamma(a,z)=z^{a-\fr}e^{-z}\bl\{d_0(-\chi) \sum_{k=0}^m \frac{A_k(\chi)}{ z^{k/2}}+\sum_{k=1}^m \frac{B_k(\chi)}{ z^{k/2}} +O(z^{-(m+1)/2})\br\}
\ee
as $z\ra\infty$ in the sector $|\arg\,z|<\fs\pi$. The function $d_0(\chi)$ is given in (\ref{e26a}) and the coefficients
$A_k(\chi)$, $B_k(\chi)$ are defined in (\ref{e29}).
\end{theorem}

The expansions (\ref{e210}) and (\ref{e211}) have been derived for $z\ra\infty$ in the sector $|\arg\,z|<\fs\pi$ uniformly in the parameter $a$. The domain of this parameter is controlled by the different restrictions on $\Re (z-a)$. These do not present any computational difficulty: if $\Re (z-a)>0$ one computes $\g(a,z)$ whereas if $\Re (z-a)<0$ one computes $\gamma(a,z)$, with the other function being determined from the identity (\ref{e11a}). If $\Re (z-a)=0$, either function can be computed. In this way we see that all values of $a$ satisfying $|\arg\,a|\leq\pi$ can be covered by one expansion or the other.

The expansions obtained contain a complementary error function of argument proportional to $\chi=(z-a)/\sqrt{z}$ that measures transition through the point $z=a$, together with two asymptotic series in inverse powers of $z^{1/2}$ with coefficients depending on $\chi$. The sectors of validity $|\arg\,z|<\fs\pi$ and $\Re (z-a)\geq 0$ or $\Re (z-a)\leq 0$ correspond to the sectors $|\arg\,(\pm\chi)|<\f{3}{4}\pi$, respectively. The form of the coefficients $A_k(\chi)$ and $B_k(\chi)$ is particularly well suited for computation in the neighbourhood of the transition point, where $|\chi|$ is not too large. For large values of $\chi$, however, the coefficients in the form (\ref{e29}) become increasingly difficult to compute on account of the severe cancellation that takes place. To illustrate, we consider the evaluation of $d_k(\chi)$ by means of the recurrence (\ref{e27}) in the case $k=4$. From (\ref{e27a}) we have
\[d_4(\chi)=(\f{1}{8}+\f{1}{4}\chi^2+\f{1}{24}\chi^4)d_0(\chi)-(\f{5}{24}\chi+\f{1}{24}\chi^3).\]
When $\chi=10$, for example, we find
\[d_4(10)=43.750008682907\ldots -43.75 \doteq 8.682907\times 10^{-6}.\]
There is an even more extreme cancellation arising in the coefficient of $z^{-2}$ in the expansion (\ref{e210})
given by $A_4(10) d_0(10)-B_4(10)$. The values of $A_4(10)d_0(10)$ and $B_4(10)$ are about $4.9637\times 10^6$, but their difference is of order $10^{-8}$. 

One way of circumventing this problem, and so extending the usefulness of the expansions (\ref{e210}) and (\ref{e211}), would be to determine the coefficients $C_k(\chi)$ from the first sum in (\ref{e28}) by computing 
the functions $d_{k+2j}(\chi)$ through use of the routine in {\em Mathematica} for the evaluation of the parabolic cylinder function. An obvious criticism of such an approach is that the parabolic cylinder function of large argument (and possibly also large order) is being computed to determine the coefficients in the expansion of the simpler incomplete gamma functions. 

A different type of uniform expansion for $a\ra\infty$ in $|\arg\,a|<\fs\pi$ has been given by Dingle \cite[p.~249]{D} who employed the same factorisation procedure applied to the integral representations in (\ref{e11}). Expressed in our notation this takes the form:
\begin{theorem}$\!\!\!.$ Let $\xi=(z-a)/\sqrt{a}$. Define the coefficients ${\hat A}_k(\xi)$ and ${\hat B}_k(\xi)$ as in (\ref{e29}) with $S_3(k,j)$ replaced by ${\hat\alpha}_j(k)$, where the ${\hat\alpha}_j(k)$ are the coefficients in the expansion of 
\[{\hat c}_k(a)=\frac{d^k}{dt^k}\,(1+t)^a e^{a(-t+\fr t^2)}|_{t=0}=\sum_{j=0}^{\lfloor k/3\rfloor} {\hat\alpha}_j(k) a^j\quad (k=0, 1, 2, \ldots).\]
Then we have 
\[\g(a+1,z)\sim z^{a+\fr}e^{-z}\bl\{d_0(\xi) \sum_{k=0}^\infty \frac{{\hat A}_k(\xi)}{ a^{k/2}}-\sum_{k=1}^\infty \frac{{\hat B}_k(\xi)}{ a^{k/2}}\br\},\quad \Re (z-a)\geq 0,\]
\[\gamma(a+1,z)\sim z^{a+\fr}e^{-z}\bl\{d_0(-\xi) \sum_{k=0}^\infty \frac{{\hat A}_k(\xi)}{ a^{k/2}}+\sum_{k=1}^\infty \frac{{\hat B}_k(\xi)}{ a^{k/2}}\br\},\quad \Re (z-a)\leq 0\]
as $a\ra\infty$ in $|\arg\,a|<\fs\pi$, where ${\hat A}_0(\xi)=1$, ${\hat A}_1(\xi)=-\xi-\f{1}{3}\xi^2$ and ${\hat B}_1(\xi)=\f{2}{3}+\f{1}{3}\xi^2$.
\end{theorem}
\vspace{0.2cm}

\noindent The coefficients ${\hat A}_k(\xi)$ and ${\hat B}_k(\xi)$ are tabulated in \cite[Table 3]{P} for $0\leq k\leq 6$.
We remark that the headings in this table are {\em correct\/} and that in the present paper we have reversed the sign of ${\hat B}_k(\xi)$ for convenience.

The nature of the expansions in Theorems 1 and 2 for large $\chi$ (that is when $|z-a|\gg |z|^{1/2}$) can be examined
by use of the result \cite[p.~347]{WW}
\[d_k(\chi)=\chi^{-k-1}\{1-\frac{(k+1)(k+2)}{2\chi^2}+O(\chi^{-4})\}\qquad (\chi\ra\infty,\ \ |\arg\,\chi|<\f{3}{4}\pi).\]
It is shown in \cite[Section 4]{P} that these expansions reduce to the forms given by Mahler \cite{M} in this limit. Similarly,
when $\xi$ is large (that is when $|z-a|\gg |a|^{1/2}$) the expansions in Theorem 3 reduce to the expansions obtained by Tricomi \cite{T}.  
\vspace{0.6cm}

\begin{center}
{\bf 3. \  The analysis of the function $Q(a,a)$}
\end{center}
\setcounter{section}{3}
\setcounter{equation}{0}
\renewcommand{\theequation}{\arabic{section}.\arabic{equation}}
In \cite{P}, the expansions in Theorems 1 and 2 were shown to satisfy the identity (\ref{e11a}) when $z=a$ ($\chi=0$).
Here we examine in more detail the expansion of $\g(a,a)$ for large $a$ and demonstrate that (\ref{e210}) reduces to
the known expansion cited in \cite[Eq. (8.12.15)]{DLMF}.

When $\chi=0$ we have 
\[d_k(0)=\frac{\sqrt{\pi}}{2^{(k+1)/2} \g(\fs k+1)}\]
so that $d_0(\chi)=\sqrt{\pi/2}$ and
\[p_{2k}(0)=\frac{1}{2^k k!},\qquad q_{2k+1}(0)=\frac{\sqrt{\pi}}{2^{k+1} \g(k+\f{3}{2})}.\]
Since $A_{2k+1}(0)=B_{2k}(0)=0$, we then find, from (\ref{e210}),
\[Q(a,a)= \frac{a^{a-\fr}e^{-a}}{\g(a)}\bl\{\sqrt{\frac{\pi}{2}} \sum_{k=0}^m \frac{A_{2k}(0)}{a^{k}}-\sum_{k=0}^m \frac{B_{2k+1}(0)}{ a^{k+\fr}}+O(a^{-m-1})\br\}\]
\[=\frac{\sqrt{2\pi}a^{a-\fr}e^{-a}}{\g(a)} \sum_{k=0}^m \frac{A_{2k}(0)}{a^{k}}\bl\{\frac{1}{2}-\frac{1}{\sqrt{2\pi a}}\,\sum_{k=0}^m \frac{E_k}{ a^{k}}+O(a^{-m-1})\br\}\]
as $a\ra\infty$ in $|\arg\,a|<\fs\pi$, where we have put
\[\frac{\sum_{k=0}^m B_{2k+1}(0) a^{-k}}{\sum_{k=0}^m A_{2k}(0)a^{-k}}=\sum_{k=0}^m \frac{E_k}{ a^{k}}+O(a^{-m-1}).\]
From (\ref{e29}), the coefficients $A_{2k}(0)$ and $B_{2k+1}(0)$ are given by
\[A_{2k}(0)=\sum_{j=0}^{2k} \frac{(-1)^j S_3(2k+2j,j)}{2^{k+j} (k+j)!},\qquad B_{2k+1}(0)=\sqrt{\pi} \sum_{j=0}^{2k+1} \frac{(-)^j S_3(2k+2j+1,j)}{2^{k+j+1} \g(k+j+\fs)}.\]

The well-known expansion of $\g(a)$ is \cite[Eq. (5.11.10), (5.11.11)]{DLMF}
\bee\label{e30}
\g(a)=\sqrt{2\pi} a^{a-\fr}e^{-a} \sum_{k=0}^m (-)^k \gamma_k a^{-k}+O(a^{-m-1})
\ee
as $a\ra\infty$ in $|\arg\,a|\leq\fs\pi$, where the first few Stirling coefficients $\gamma_k$ are
\[\gamma_0=1,\quad \gamma_1=-\f{1}{12}, \quad \gamma_2=\f{1}{288},\quad \gamma_3=\f{139}{51840},\quad \gamma_4=-\f{571}{2488320}, \ldots\ .\]
It has been shown in \cite[Theorem 2.4]{BM} that $A_{2k}(0)=(-)^k \gamma_k$. 
Consequently we have the expansion
\bee\label{e31}
Q(a,a)=\frac{1}{2}-\frac{1}{\sqrt{2\pi a}}\sum_{k=0}^m E_k a^{-k}+O(a^{-m-1})
\ee
as $a\ra\infty$ in $|\arg\,a|<\fs\pi$, where routine evaluation of the coefficients $E_k$ with the help of {\it Mathematica} shows that
\[E_0=\f{1}{3},\quad E_1=\f{1}{540} \quad E_2=-\f{25}{6048},\quad E_3=-\f{101}{155520}, \quad E_4=\f{3184811}{3695155200}, \]
\[E_5=\f{2745493}{8151736320}, \quad E_6=-\f{119937661}{225740390400},\quad E_7=-\f{8325705316049}{24176795811840000},\ldots .\]
The expansion (\ref{e31}) agrees with that given in \cite[Eq. (8.12.15), (8.12.17)]{DLMF}, where values of $E_k$ are given up to $k=5$.
\vspace{0.6cm}

\begin{center}
{\bf 4. \  Concluding remarks}
\end{center}
\setcounter{section}{4}
\setcounter{equation}{0}
\renewcommand{\theequation}{\arabic{section}.\arabic{equation}}
We have shown that the table of the coefficients $A_k(\chi)$ and $B_k(\chi)$ that appear in the asymptotic expansions of the incomplete gamma functions given in \cite[Table 2]{P} has incorrect headings, although the forms of these coefficients themselves are correct.  This has led to an error in the presentation of these expansions in Eqs. (8.12.18) -- (8.12.20) in the DLMF Handbook \cite{DLMF}. This entry should read:
\[\left.\begin{array}{l}Q(a,z)\\P(a,z)\end{array}\!\!\right\} \sim \frac{z^{a-\fr}e^{-z}}{\g(a)}\bl\{d(\pm\chi)\sum_{k=0}^\infty \frac{A_k(\chi)}{z^{k/2}}\mp \sum_{k=1}^\infty \frac{B_k(\chi)}{z^{k/2}}\br\}\]
as $z\ra\infty$ in $|\arg\,z|<\fs\pi$, and $\Re (z-a)\leq 0$ for $P(a,z)$ and $\Re (z-a)\geq 0$ for $Q(a,z)$.
Here
\[\chi=\frac{z-a}{\sqrt{z}},\qquad d(\pm\chi)=\sqrt{\frac{\pi}{2}}\,e^{\chi^2/2} \,\mbox{erf} \bl(\frac{\pm\chi}{\sqrt{2}}\br)\]
and $A_0(\chi)=1$, $A_1(\chi)=\fs\chi+\f{1}{6}\chi^3$, $B_1(\chi)=\f{1}{3}+\f{1}{6}\chi^2$.

We have also examined more closely the expansion for $Q(a,a)$ and have shown that
(\ref{e210}) reduces to the expansion cited in \cite[Eq. (8.12.15)]{DLMF}.

\end{document}